\input amstex
\input epsf
\magnification=\magstep1 
\baselineskip=13pt
\documentstyle{amsppt}
\vsize=8.7truein \CenteredTagsOnSplits \NoRunningHeads
\def\EE{\bold {E}}
\def\dist{\operatorname{dist}}

\def\ii{\bold {i}}
\topmatter
 
\title On the zeros of partition functions with multi-spin interactions \endtitle 
\author Alexander Barvinok  \endauthor
\address Department of Mathematics, University of Michigan, Ann Arbor,
MI 48109-1043, USA \endaddress
\email barvinok$\@$umich.edu \endemail
\date June 27, 2024\enddate
\thanks  This research was partially supported by NSF Grant DMS 1855428. 
\endthanks 
\abstract Let $X_1, \ldots, X_n$ be probability spaces, let $X$ be their direct product, let $\phi_1, \ldots, \phi_m: X \longrightarrow {\Bbb C}$ be random 
variables, each depending only on a few coordinates of a point $x=(x_1, \ldots, x_n)$, and let $f=\phi_1 + \ldots + \phi_m$. The expectation $\EE\thinspace e^{\lambda f}$, where $\lambda \in {\Bbb C}$, appears in statistical physics as the partition function of a system with multi-spin interactions, and also in combinatorics and computer science, where it is known as the partition function of edge-coloring models, tensor network contractions or a Holant polynomial. Assuming that each $\phi_i$ is 1-Lipschitz in the Hamming metric of $X$, that each $\phi_i(x)$ depends on at most $r \geq 2$ coordinates $x_1, \ldots, x_n$ of $x \in X$, and that for each $j$ there are at most $c \geq 1$ functions $\phi_i$ that depend on the coordinate $x_j$, we prove that $\EE\thinspace e^{\lambda f} \ne 0$ provided $| \lambda | \leq \ (3 c \sqrt{r-1})^{-1}$ and that the bound is sharp up to a constant factor. Taking a scaling limit, we prove a similar result for functions $\phi_1, \ldots, \phi_m: {\Bbb R}^n \longrightarrow {\Bbb C}$ that are 1-Lipschitz in the $\ell^1$ metric of ${\Bbb R}^n$ and where the expectation is taken with respect to the standard Gaussian measure in ${\Bbb R}^n$.
As a corollary, the value of the expectation can be efficiently approximated, provided $\lambda$ lies in a slightly smaller disc. 
\endabstract
\endtopmatter
\document

\head 1. Introduction and the main results \endhead

We investigate the complex zeros and computational complexity of functionals of a particular type, which appear under different names in statistical physics, computer science, and combinatorics.  The functionals can be viewed as the partition function of a spin system with multiple spin interactions (statistical physics) or as the partition function of a hypergraph edge-coloring model, also known as a tensor network contraction, or as a Holant polynomial (combinatorics and computer science).

In what follows, we consider functions on the direct product $X_1 \times \ldots \times X_n$ of probability spaces. For a point $x \in X$, $x=\left(x_1, \ldots, x_n\right)$, we refer to $x_j \in X_j$ as the $j$-th {\it coordinate} of $x$. The {\it Hamming distance} between two points $x, y \in X$ is the number of the coordinates where they differ:
$$\dist(x, y)=\left|j:\ x_j \ne y_j \right|, \quad \text{where} \quad x=\left(x_1, \ldots, x_n \right) \quad \text{and} \quad y=\left(y_1, \ldots, y_n \right).$$

We prove the following main result.
\proclaim{(1.1) Theorem} Let $X_1, \ldots, X_n$ be probability spaces, let $X=X_1 \times \ldots \times X_n$ be the product space and let 
$\phi_1, \ldots, \phi_m: X \longrightarrow {\Bbb C}$ be measurable functions. 
Assume that 
\roster
\item Each function  $\phi_i$ is 1-Lipschitz in the Hamming metric, that is, 
$$\left| \phi_i(x)-\phi_i(y)\right| \leq 1$$
provided $x=\left(x_1, \ldots, x_n\right)$ and $y=\left(y_1, \ldots, y_n\right)$ differ in one coordinate;
\item Each function $\phi_i$ depends only on at most $r \geq 2$ coordinates, that is, for each $i=1, \ldots, m$ there is a subset $J_i \subset \{1, \ldots, n\}$ with 
$|J_i| \leq r$ such that 
$$\phi_i\left(x_1, \ldots, x_n\right)=\phi_i\left(y_1, \ldots, y_n \right) \quad \text{whenever} \quad x_j =y_j \quad \text{for all} \quad j \in J_i;$$
\item For every $j=1, \ldots, n$, there are at most $c \geq 1$ functions $\phi_i$ that depend on the coordinate $x_j$, that is, 
$\left|i: \ j \in J_i \right| \leq c$ for all $j$.
\endroster
Let 
$$f=\sum_{i=1}^m \phi_i$$ and 
suppose that $\lambda$ is a complex number such that 
$$|\lambda| \ \leq \ {1 \over 3 c \sqrt{r-1}}.$$
Then 
$$\EE\thinspace e^{\lambda f} \ne 0.$$
Moreover, if, additionally, 
$$\left| \phi_i(x)\right| \ \leq \ L \quad \text{for all} \quad x \in X \quad \text{and} \quad i=1, \ldots, m \tag1.1.1$$
and some $L > 0$, then 
$$e^{|\lambda| m L} \ \geq \ \left|\EE\thinspace e^{\lambda f}\right| \ \geq \ e^{-|\lambda| m L} \left( \cos {\pi \over 4 \sqrt{r-1}}\right)^n. \tag1.1.2$$
\endproclaim
We prove Theorem 1.1 in Section 2, and in Section 4 we show that the bound for $|\lambda|$ is optimal, up to a constant factor. The method of proof extends those of \cite{Ba17} and \cite{BR19}. The dependence on $r$ is worth noting. One approach frequently used for problems of this type is the {\it cluster expansion} method, see \cite{Je24} for a recent survey. In the situation of Theorem 1.1 it apparently gives only $|\lambda|=\Omega(1/rc)$ as a bound for the zero-free region, while also requiring $|\phi_i|$ to remain uniformly bounded \cite{Ca+22}.

The dependence on $r$  allows us to obtain a zero-free region for a partition function which can be considered as a scaling limit of that of Theorem 1.1. We consider the standard Gaussian probability measure in ${\Bbb R}^n$ with density 
$${1 \over (2\pi)^{n/2}} \exp\left\{ -{1 \over 2} \left(x_1^2 + \ldots + x_n^2\right)\right\} \quad \text{for} \quad x=\left(x_1, \ldots, x_n\right)$$
and prove the following result.

\proclaim{(1.2) Theorem} Let $\phi_1, \ldots, \phi_m: {\Bbb R}^n \longrightarrow {\Bbb C}$ be functions on Euclidean space ${\Bbb R}^n$, endowed with 
the standard Gaussian probability measure. Assume that 
\roster 
\item Each function $\phi_i$ is 1-Lipschitz in the $\ell^1$ metric of ${\Bbb R}^n$, that is, 
$$\left| \phi_i\left(x_1, \ldots, x_n\right)-\phi_i\left(y_1, \ldots, y_n \right)\right| \ \leq \ \sum_{i=1}^n \left| x_i-y_i \right|;$$
\item Each function $\phi_i$ depends only on at most $r \geq 2$ coordinates;
\item For every $j=1, \ldots, n$, there are at most $c \geq 1$ functions $\phi_i$ that depend on the coordinate $x_j$;
\endroster
Let 
$$f=\sum_{i=1}^m \phi_i$$ and 
suppose that $\lambda$ is a complex number such that 
$$|\lambda| \ < \ {1 \over 6 c \sqrt{r-1}}.$$
Then 
$$\EE\thinspace e^{\lambda f} \ne 0.$$
Moreover, if, additionally, 
$$\left| \phi_i(x)\right| \ \leq \ L \quad \text{for all} \quad x \in {\Bbb R}^n \quad \text{and} \quad i=1, \ldots, m \tag1.2.1$$
and some $L > 0$, then 
$$e^{|\lambda| m L} \ \geq \ \left|\EE\thinspace e^{\lambda f}\right| \ \geq \ e^{-|\lambda| m L} e^{-\pi^2 n/32 r}. \tag1.2.2$$
\endproclaim

We prove Theorem 1.2 in Section 3.

Suppose we want to approximate the value of $\EE\thinspace e^{\lambda f}$ efficiently. As we are dealing with complex numbers, it is convenient to adopt the following definition: we say that a complex number $z_1 \ne 0$ approximates another complex number $z_2 \ne 0$ within relative error $\epsilon$ if 
$z_1=e^{w_1}$ and $z_2=e^{w_2}$ for some $w_1$ and $w_2$ such that $|w_1-w_2| \leq \epsilon$.

Let us fix some $0 < \delta < 1$. From by now standard method of polynomial interpolation \cite{Ba16}, \cite{PR17}, it follows that for any $\lambda$
such that 
$$|\lambda| \ \leq \ {1-\delta \over 3 c \sqrt{r-1}} \quad \text{(in Theorem 1.1)} \quad \text{or} \quad |\lambda| \ \leq \ {1-\delta \over 6 c \sqrt{r-1}} \quad
\text{(in Theorem 1.2)}$$ 
approximating the value of $\EE\thinspace e^{\lambda f}$ within relative error $0 < \epsilon < 1$ reduces to the computation of the $m^k$ expectations
$$\EE \left(\phi_{i_1} \cdots \phi_{i_k}\right) \quad \text{for all} \quad 1 \leq i_1, i_2, \ldots, i_k \leq m.$$
where 
$$k=O_{\delta}\left(\ln(m+n)-\ln \epsilon\right).$$ We sketch the algorithm in Section 5. In many applications of Theorem 1.1 the probability spaces $X_j$ are finite. Assuming that each $X_j$ contains at most $q$ elements, computing each expectation $\EE \left(\phi_{i_1} \cdots \phi_{i_k}\right)$
by the direct enumeration takes $O\left(q^{rk}\right)$ time, provided we have an access to the values of $\phi_i(x)$ for a given $x \in X$. If the parameter $r$ is fixed in advance, we obtain a quasi-polynomial algorithm of $(q(m+n))^{O(\ln(m+n) -\ln \epsilon)}$ complexity to approximate $\EE\thinspace e^{\lambda f}$ within relative error $0< \epsilon < 1$. Moreover, the general technique of Patel and Regts \cite{PR17} allows one to speed it up to a genuinely polynomial time algorithm of $\bigl(q(m+n)\epsilon^{-1}\bigr)^{O(1)}$ complexity, provided the parameter $c$ is also fixed in advance. In the context of Theorem 1.2, the expectations $\EE \left(\phi_{i_1} \cdots \phi_{i_k}\right)$ are represented by integrals in the space of dimension at most $kr$ and often can be efficiently computed or approximated with high accuracy, assuming again that $r$ is fixed in advance. In that case, we also obtain an algorithm of quasi-polynomial complexity to approximate $\EE\thinspace e^{\lambda f}$.

As we mentioned, the expectation $\EE\thinspace e^{\lambda f}$ appears in several different, though closely related contexts. 

\subhead (1.3) Statistical physics: partition functions of systems with multi-spin interactions \endsubhead Here we describe the statistical physics context of Theorem 1.1.
Suppose that we have a system of $n$ particles, where the $j$-th particle can be in a state described by a point $x_j \in X_j$, also called the {\it spin} of the particle. The vector 
$x=\left(x_1, \ldots, x_n\right)$ of all spins is called a {\it spin configuration}. The particles interact in various ways, and the energy of a spin configuration is given by a function 
$-f(x_1, \ldots, x_n)$. Then for a real $\lambda >0$, interpreted as the {\it inverse temperature}, the value of $\EE\thinspace e^{\lambda f}$ is the {\it partition function} of the system, see, for example,  \cite{FV18} for the general setup. If the function $f$ is written as a sum of $\phi_i$, as in Theorem 1.1, then the energy of the system is a sum over subsystems, each containing at most $r$ particles, and each particle participating in at most $c$ subsystems.

The idea to relate complex zeros of the partition function $\lambda \longmapsto \EE\thinspace e^{\lambda f}$ to the physical phenomenon of {\it phase transition} apparently goes back to Mayer \cite{Ma37}, when it became known as ``Mayer's conjecture", see Section 4.12.3 of \cite{FV18} for a discussion.
It then took off with the classical works of Lee and Yang \cite{YL52}, \cite{LY52}, who showed that zero-free regions like the one described by Theorem 1.1 correspond to regimes with no phase transition. 

In terms of Theorem 1.1, the most studied case is that of $r=2$, which includes the classical Ising, Potts and Heisenberg models, see \cite{FV18}. In that case, we have a (finite) graph $G=(V, E)$ with set $V$ of vertices and set $E$ of edges. The particles are identified with the vertices of $G$,  the spins are $\pm 1$ in the case of the Ising model, elements of some finite set in the case of the Potts model, or vectors in Euclidean space, in the case of the Heisenberg model. The interactions are pairwise and described by the functions attached to the edges of $E$. Hence, in the context of Theorem 1.1, we have $r=2$ and $c$ is the largest degree of a vertex of $G$. Starting with \cite{LY52}, the complex zeros of the partition function of systems with pairwise interactions were actively studied in great many papers. We refer to \cite{FV18} for earlier, and to \cite{B+21}, \cite{G+22}, \cite{PR20}, \cite{P+23} for more recent works.

Zero-free regions for systems with multiple spin interactions were considered by Suzuki and Fisher \cite{SF71}, again in connection with phase transition, see also \cite{L+19} and \cite{L+16} for recent developments. In this case, the particles are identified with the vertices of a hypergraph, while functions attached to the edges (sometimes called hyperedges) of the hypergraph describe interactions. In terms of Theorem 1.1, the maximum number of vertices of an edge of the hypergraph is $r$, while $c$ is the largest degree of vertex. The papers \cite{SF71} and \cite{L+19}, respectively \cite{L+16}, consider rather special {\it ferromagnetic}, respectively {\it antiferromagnetic}, types of interaction, so their results are not directly comparable to ours. While \cite{SF71}, \cite{L+19} and \cite{L+16} say more about those specific models, our Theorem 1.1 allows one to handle a wider class of interactions.

\subhead (1.4) Statistical physics: partition functions of interacting particles with an external field  \endsubhead
Here we give an example of how Theorem 1.2 applies to systems of pairwise interacting particles, cf. \cite{FV18}. For an integer $N > 1$, let us consider $N$ particles represented as
vectors $x^{(1)}, \ldots, x^{(N)}$ in ${\Bbb R}^d$, chosen independently at random from the standard Gaussian distribution in ${\Bbb R}^d$. Suppose that there are pairwise repulsive constant forces, so that the total energy of the system is 
 $$-f\left(x^{(1)}, \ldots, x^{(N)}\right)= -\sum_{i \ne j} \|x^{(i)}-x^{(j)}\|,$$
 where $\| \cdot \|$ is the standard Euclidean norm in ${\Bbb R}^d$.
 The energy is minimized when the particles are far away from each other. However, the probability of a configuration with large pairwise distances is small, so that the Gaussian density can be interpreted as an external force pushing the particles towards the origin,
 In terms of Theorem 1.2, we have 
$n=dN$, $m={N \choose 2}$, $r=2d$ and $c=N$. Theorem 1.2 establishes a zero-free region of the order
$$|\lambda| = \Omega\left({1 \over N \sqrt{d}}\right)$$
for the partition function $\EE\thinspace e^{\lambda f}$.
This model is related to the old problem of finding a configuration of points on the unit sphere in ${\Bbb R}^d$ that maximizes the sum of pairwise distances between points, see \cite{B+23} (we note that for large $d$ the standard Gaussian measure in ${\Bbb R}^d$ is concentrated around the sphere $\|x\|=\sqrt{d}$).
\bigskip
Essentially identical models described by Theorem 1.1 are studied in connection with questions in combinatorics and theoretical computer science.

\subhead (1.5) Combinatorics and computer science: edge-coloring models and tensor network contractions \endsubhead
Let $G=(V, E)$ be a graph. Suppose now that every edge of $G$ can be in one of the $k$ states $\{1, \ldots, k\}$, typically called ``colors". In terms of 
Theorem 1.1, we have $n=|E|$ and $X_1=\ldots = X_n = \{1, \ldots, k\}$, so that $X =\{1, \ldots, k\}^E$. Suppose further, that to each vertex $v$ of $G$ a function 
$\psi_v: X_1 \times \ldots \times X_ n \longrightarrow {\Bbb C}$ is attached, that depends only on the colors of the edges of $G$ that contain $v$.
The expression
$$\sum_{x \in X} \prod_{v \in V} \psi_v(x) \tag1.5.1$$
is known as the {\it tensor network} contraction, or the partition function of the {\it edge coloring model}, or as a {\it Holant polynomial}, see \cite{Re18} for relations between different models, and references. For relations with spin systems in statistical physics and also knot invariants, see \cite{HJ93}. Assuming that 
$\psi_v(x) \ne 0$ for all $x$ and $v$, we can write $\psi_v(x)=\exp\left\{ \phi_v(x)\right\}$ and (1.5.1) can be written as the scaled expectation $\EE\thinspace e^{\lambda f}$ of Theorem 1.1, assuming the uniform probability measure in each space $X_j$.
In the context of Theorem 1.1, we have $c=2$, while $r$ is equal to the largest degree of a vertex of $G$. Note that compared to the graph interpretation of the Ising and Potts models from Section 1.3, the parameters $r$ and $c$ switch places.

A similar to (1.5.1) expression can be built for hypergraphs, in which case $r$ is still the maximum degree of vertex, while $c$ becomes the maximum size of an edge, see for example, Chapter 6 of \cite{Ta11} for the discussion of the partition function of a $k$-satisfiability problem.

In combinatorics and computer science, there is a lot of interest in efficient approximation of (1.5.1), also in connection with zero-free regions \cite{GG16}, \cite{B+22}, \cite{Ca+22}, \cite{G+21}, \cite{Re18}, since many interesting counting problems on graphs and hypergraphs can be stated as a problem of computing (1.5.1). To illustrate, we consider just one example of perfect matchings in a hypergraph, see \cite{Ta11} for further motivation.

Let $H=(V, E)$ be a hypergraph with set $V$ of vertices and set $E$ of edges. We choose $X=\{1, 2\}^E$ so that every edge of $H$ can be colored into one of the two colors, which we interpret as the edge being selected or not selected. We define $\psi_v(x)=1$ if exactly one edge containing $v$ is selected, and 
$\psi_v(x)=0$ otherwise. Then (1.5.1) is exactly the number of {\it perfect matchings} is $H$. Since deciding whether a hypergraph contains a perfect matching is a well-known NP-complete problem, there is little hope to approximate (1.5.1) efficiently. One can try to come up with a more approachable version of the problem by modifying the definition of $\psi_v$ so that
$\psi_v(x)=1$ if exactly one edge containing $v$ is selected, and $\psi_v(x)=1-\delta$ otherwise, for some $0 < \delta < 1$. In this case, the sum (1.5.1), while taken over all collections of edges of $H$, is  ``exponentially tilted" towards perfect matchings. Thus every perfect matching is counted with weight 1,
and an arbitrary collection of edges is discounted exponentially in the number of vertices where the perfect matching condition is violated.
In other words, the weight  of a collection of edges is exponentially small in the number of vertices that belong to any number of edges in the collection other than 1.

The larger $\delta$ we are able to choose, the more (1.5.1) is tilted towards perfect matchings. It follows from Theorem 1.1 via the interpolation method that 
for a fixed $r$, there is a quasi-polynomial  algorithm approximating (1.5.1) for some 
$\delta=\Omega\left(1/c \sqrt{r}\right)$,
where $c$ is the maximum number of vertices in an edge of $H$ and $r$ is the largest degree of vertex. The results \cite{Ca+22} and \cite{Re18} appear to allow only for
$\delta=\Omega\left(1/cr\right)$.
 We note that Theorem 1.1 allows us to choose 
$$\psi_v=\exp\left\{ -\Omega\left(|k-1|\over c \sqrt{r}\right)\right\}, $$ where $k$ is the number of selected edges containing $v$, so as to assign smaller weights to the vertices $v$ with bigger violation of the perfect matching condition, up to the smallest weight of $\psi_v =\exp\left\{ -\Omega(\sqrt{r}/c)\right\}$, when all edges containing $v$ are selected.

We note also that Theorem 1.1 allows us to select edges with a non-uniform probability, and hence to zoom in the perfect matchings even more: if the hypergraph $H$ is $r$-regular, that is, if each vertex is contained in exactly $r$ edges, it makes sense to select each edge independently with probability $1/r$, so that for each vertex $v$ the expected number of selected edges containing $v$ is exactly 1. If $H$ is not regular, one can choose instead the {\it maximum entropy distribution}, which also ensures that the expected number of selected edges containing any given vertex is exactly 1, while the weights of all perfect matchings remain equal.  The maximum entropy distribution exists if an only if there exists a positive {\it fractional perfect matching}, that is, an assignment of positive real weights to the edges of $H$ such that for every vertex of $H$ the sum of weights of the edges containing it is exactly 1, see \cite{Ba23} for details.

\subhead (1.6) Numerical integration in higher dimensions \endsubhead Computationally efficient numerical integration of multivariate functions  is an old problem that is often associated with ``the curse of dimensionality". The most spectacular success is achieved for integration of log-concave functions 
on ${\Bbb R}^n$ via the Monte Carlo Markov Chain method, see \cite{LV07} for a survey. Deterministic methods enjoyed much less success. In \cite{GS24}, the authors, using the deterministic {\it decay of correlations} approach, considered integration of functions over the unit cube $[0, 1]^n$. The model of \cite{GS24} fits the setup of our Theorem 1.1, if we choose $X_1= \ldots = X_n =[0,1]$ endowed with the Lebesgue probability measure. The results of \cite{GS24} appear to be weaker than our Theorem 1.1, in the dependence on both parameters $r$ and $c$, as well as in the class of allowed functions $\phi_i$.

Theorem 1.2 allows us to integrate efficiently some functions that are decidedly not log-concave (or log-convex), for example if we choose 
$\phi_i(x)=\bigl| |x_i|-1 \bigr|$ for some $i$, thus reaching beyond the realm of functions efficiently integrated by randomized methods.

\subhead (1.7) Other applications \endsubhead Zero-free regions of partition functions of the type covered by Theorem 1.1 turn out to be relevant to the {decay of correlations} \cite{Ga23}, \cite{Re23}, to the mixing time of Markov Chains \cite{Ch+22}, to the validity of the Central Limit Theorem for 
combinatorial structures \cite{MS19}, as well as to other related algorithmic applications \cite{J+22}.

\head 2. Proof of Theorem 1.1 \endhead 

We start with a lemma.

\proclaim{(2.1) Lemma} Let $X$ be a probability space and let $f: X \longrightarrow {\Bbb C}$ be an integrable function. Suppose that $f(x) \ne 0$ for all 
$x \in X$ and, moreover, for any two points $x, y \in X$, the angle between $f(x)$ and $f(y)$, considered as non-zero vectors in ${\Bbb R}^2={\Bbb C}$ does not exceed $\theta$ for some $0 \leq \theta < 2\pi/3$. Then 
$$\left| \EE\thinspace f \right| \ \geq \ \left(\cos {\theta \over 2} \right) \EE\thinspace |f|.$$
\endproclaim
\demo{Proof} This is Lemma 3.3 from \cite{BR19}, see also Lemma 3.6.3 of \cite{Ba16}. The idea of the proof is to compare both sides of the inequality to the expected length of the orthogonal projection of $f(x)$ onto the bisector of the angle enclosing the range of $f$.
{\hfill \hfill \hfill} \qed
\enddemo

We will need some technical inequalities.

\proclaim{(2.2) Lemma} The following inequalities hold
\roster
\item 
$$\left( \cos {\alpha \over \sqrt{k}}\right)^k \ \geq \ \cos \alpha \quad \text{for} \quad 0 \leq \alpha \leq {\pi \over 2} \quad \text{and} \quad k \geq 1;$$
\item
$$\sin (\tau \alpha) \ \geq \ \tau \sin \alpha \quad \text{for} \quad 0 \leq \alpha \leq {\pi \over 2} \quad \text{and} \quad 0 \leq \tau \leq 1;$$
\item 
$$\left| e^z-1 \right| \ \leq \ 2|z| \quad \text{for} \quad z \in {\Bbb C} \quad \text{such that} \quad |z| \leq 1.$$
\endroster
\endproclaim
\demo{Proof} The inequality of Part (1) is proved, for example, in \cite{Ba23}, see Lemma 6.1 there. To prove the inequality of Part (2), let 
$$f(\alpha)=\sin (\tau \alpha) \quad \text{and} \quad g(\alpha)=\tau \sin \alpha.$$
Then $f(0)=g(0)=0$ and 
$$f'(\alpha)=\tau \cos (\tau \alpha) \ \geq \ \tau \cos \alpha = g'(\alpha),$$
from which the proof follows.
To prove the inequality of Part (3), we note that 
$$\left| e^{z}-1 \right| =\left| \sum_{k=1}^{\infty} {z^k \over k!} \right| \ \leq \ \sum_{k=1}^{\infty} {|z|^k \over k!},$$
and hence it suffices to check the inequality assuming that $z$ is non-negative real. Since the function 
$${e^z-1 \over z}=\sum_{k=0}^{\infty} {z^k \over (k+1)!}$$
is increasing for $z \geq 0$, it suffices to check the inequality at $z=1$, where it states that $|e-1| \leq 2$.
{\hfill \hfill \hfill} \qed
\enddemo

Let $X_1, \ldots, X_n$ be probability spaces and let $f: X_1 \times \ldots \times X_n \longrightarrow {\Bbb C}$ be an integrable function. For a subset 
$S \subset \{1, \ldots, n\}$, by $\EE_S f$ we denote the conditional expectation obtained by integrating $f$ over the coordinates $x_j \in X_j$ with $j \in S$. Thus $\EE_S f$ is a function of the coordinates $x_j$ with $j \notin S$. In particular, for $S=\emptyset$, we have
$\EE_S f=f$ and for $S=\{1, \ldots, n\}$, we have $\EE_S =\EE\thinspace f$. 

Now we are ready to prove Theorem 1.1.

\subhead (2.3) Proof of Theorem 1.1 \endsubhead  Let
$$\theta= {\pi \over 2 \sqrt{r-1}}. \tag2.3.1$$
For a subset $S \subset \{1, \ldots, n\}$, we prove by induction on $|S|$ the following 
\bigskip
\noindent (2.3.2) {\sl Statement:} Let functions $\phi_i$, $f$ and a complex number $\lambda$ be as in Theorem 1.1.
 For every $S \subset \{1, \ldots, n\}$, we have 
$\EE_S\thinspace e^{\lambda f} \ne 0$, by which we mean that $\EE_S e^{\lambda f}  \ne 0$ for any choice of the coordinates $x_j \in X_j$ with $j \notin S$. Moreover, 
for every $j \notin S$, the value of $\EE_S\thinspace e^{\lambda f} \ne 0$, considered as a vector in ${\Bbb R}^2={\Bbb C}$, rotates by not more than an angle of  $\theta$ when only the $x_j$ coordinate of $x=\left(x_1, \ldots, x_n\right)$ changes, while all other coordinates stay the same.
\bigskip
Once we have (2.3.2) for $S=\{1, \ldots, n\}$, we conclude that $\EE\thinspace e^{\lambda f} \ne 0$, which is what we want to prove.

Suppose that $|S|=0$, so that $S=\emptyset$ and 
$$\EE_S \thinspace e^{\lambda f}= e^{\lambda f}=\exp\left\{ \lambda \sum_{i=1}^m \phi_i \right\}.$$
Let $x', x'' \in X$ be two points that differ only in the $x_j$ coordinate. Since at most $c$ of the functions $\phi_i$ depend on the coordinate $x_j$ and each function $\phi_i$ is 1-Lipschitz, we have 
$$|\lambda f(x') - \lambda f(x'')| =|\lambda| |f(x')-f(x'')| \ \leq \  |\lambda| c \ \leq \ {1 \over 3\sqrt{r-1}} \ < \ \theta,$$
which establishes the base of the induction.

Suppose now that  Statement 2.3.2 holds for every $S$ with $|S| \leq k$, where $0 \leq k < n$. Let us choose an arbitrary $S \subset \{1, \ldots, n\}$ such that 
$|S|=k+1$ and an index $j \in S$.
Applying Lemma 2.1 and the induction hypothesis, we conclude that 
$$\left| \EE_S\thinspace e^{\lambda f} \right| =\left| \EE_{\{j\}} \EE_{S\setminus \{j\}} e^{\lambda f} \right| 
\ \geq \ \left(\cos {\theta \over 2}\right) \EE_{\{j\}} \left| \EE_{S \setminus \{j\}} e^{\lambda f}\right|. \tag2.3.3$$
It follows that 
$$\EE_S\thinspace e^{\lambda f} \ne 0.$$ Assuming that $k+1 < n$, let us pick some index not in $S$, without loss of generality index $n$. We need to prove that as the coordinate 
$x_n$ changes from some value $x_n'$ to some value $x_n''$, while other coordinates remain the same, the value of $\EE_S\thinspace e^{\lambda f}$ rotates through an angle of at most $\theta$.
Let
$$I=\left\{i: \ \phi_i \quad \text{depends on} \quad x_n \right\}, \quad \text{so} \quad |I| \leq c.$$
Without loss of generality, $I \ne \emptyset$. For each $i \in I$, we define two functions 
$\phi_i'$ and $\phi_i''$, obtained from $\phi_i$ by fixing the coordinate $x_n$ to $x_n'$ and $x_n''$ respectively. Although $\phi_i'$ and $\phi_i''$ are functions of the first $n-1$ coordinates $x_1, \ldots, x_{n-1}$ of $x \in X$, we formally consider them as functions on $X$, by ignoring the last coordinate $x_n$.

Since each function $\phi_i$ is 1-Lipschitz in the Hamming metric, we have 
$$\left| \phi_i'(x)-\phi_i''(x) \right|  \ \leq \ 1 \quad \text{for all} \quad x \in X \quad \text{and} \quad i \in I. \tag2.3.4$$
Thus the value of $\EE_S\thinspace e^{\lambda f}$, where we fix $x_n=x_n'$, is 
$$\EE_S \exp\left\{ \lambda \sum_{i \in I} \phi_i' + \lambda \sum_{i \notin I} \phi_i \right\}, \tag2.3.5$$
while the value of $\EE_S e^{\lambda f}$, where we fix $x_n=x_n''$, is 
$$ \EE_S \exp\left\{ \lambda \sum_{i \in I} \phi_i'' + \lambda \sum_{i \notin I} \phi_i \right\}. \tag2.3.6$$ 
We pass from (2.3.5) to (2.3.6) step by step, replacing one $\phi_i'$ by $\phi_i''$ at each step.
Our goal is to prove that at each step, the expectation rotates by at most $\theta/c$. Once we prove that, it would follow that the expectation 
$\EE_S\thinspace e^{\lambda f}$ rotates by at most $\theta$, when we replace the value of $x_n=x_n'$ by $x_n=x_n''$. 

Let us pick an index in $I$, without loss of generality index $m$. We define 
$$f'=\phi_m' +  \psi  \quad \text{and} \quad 
f''= \phi_m'' + \psi, \tag2.3.7$$
where $\psi$ is the sum of some functions $\phi_i'$, $\phi_i''$ and $\phi_i$, where for each $i \ne m$ we select exactly one of the three functions 
$\phi_i'$, $\phi_i''$ or $\phi_i$ into the sum for $\psi$. Hence our goal is to show that the angle between 
$\EE_S\thinspace e^{\lambda f'} $ and $\EE_S\thinspace e^{\lambda f''}$ does not exceed $\theta/c$.

Since each of the functions $\phi_i'$ and $\phi_i''$ is obtained from $\phi_i$ by specifying some value of $x_n$, we can apply the induction hypothesis both to 
$f'$ and to $f''$. Let $S_0 \subset S$ be the set of indices $j \in S$ such that $\phi_m'$ and $\phi_m''$ depend on $x_j$. In particular, $|S_0| \leq r-1$. 

If $S_0 = \emptyset$, then 
$$\EE_S\thinspace e^{\lambda f'} = e^{\lambda \phi_m'} \EE_S\thinspace e^{\lambda \psi} \quad \text{and} \quad \EE_S\thinspace e^{\lambda f''} =e^{\lambda \phi_m''}
\EE_S\thinspace e^{\lambda \psi}.$$
Using  (2.3.4), we conclude the angle between $\EE_S\thinspace e^{\lambda f'} \ne 0 $ and $\EE_S\thinspace e^{\lambda f''} \ne 0$ does not exceed
$$|\lambda \phi_m''(x) - \lambda \phi_m' (x)| \ \leq \ |\lambda| \ \leq \ {1 \over  3c \sqrt{r-1}} \ < \ {\theta \over c}.$$

Suppose now that $S_0 \ne \emptyset$. 
Since $\phi_m'$ and $\phi_m''$ do not depend on the coordinates $x_j$ with $j \notin S_0$, from (2.3.7) we have
$$\split\EE_S\thinspace e^{\lambda f''}=&\EE_S \left(e^{\lambda \phi_m''} e^{\lambda \psi}\right)=\EE_S \left(e^{ \lambda \phi_m'' -\lambda \phi_m' }
e^{\lambda \phi_m'+ \lambda \psi}\right)=\EE\left(e^{ \lambda \phi_m'' -\lambda \phi_m'} e^{\lambda f'} \right) \\=
&\EE_{S_0} \left(\left( e^{\lambda \phi_m'' - \lambda \phi_m'} \right) \EE_{S\setminus S_0} e^{\lambda f'}\right) \endsplit$$
and hence 
$$\split  \EE_S\thinspace e^{\lambda f''} - \EE_S\thinspace e^{\lambda f'} =
&\EE_{S_0} \left(\left( e^{\lambda \phi_m'' - \lambda \phi_m'} \right) \EE_{S\setminus S_0} e^{\lambda f'}\right)- \EE_{S_0}\left( \EE_{S\setminus S_0} e^{\lambda f'}\right)\\ =
 &\EE_{S_0} \left( \left( e^{\lambda \phi_m''-\lambda \phi_m'}-1 \right) \EE_{S\setminus S_0}\thinspace
e^{\lambda f'} \right). \endsplit$$
By (2.3.4) and
by Part (3) of Lemma 2.2, we have 
$$\left| e^{\lambda \phi_m''-\lambda \phi_m'}-1\right| \ \leq \ 2|\lambda|.$$
Therefore,
$$\left| \EE_S\thinspace e^{\lambda  f''} - \EE_S\thinspace e^{\lambda f'} \right| \ \leq \  2|\lambda| \EE_{S_0} \left| \EE_{S \setminus S_0} \thinspace e^{\lambda f'} \right|.  \tag2.3.8$$
Iterating (2.3.3) with $f$ replaced by $f'$, we obtain that 
$$\aligned  \left| \EE_S e^{\lambda f'} \right| = &\left| \EE_{S_0} \EE_{S\setminus S_0}\thinspace e^{\lambda f'}\right| \ \geq \ \left( \cos {\theta \over 2}\right)^{|S_0|} \EE_{S_0}  \left| \EE_{S\setminus S_0} e^{\lambda f' }  \right|  \\ & \geq \ \left( \cos {\theta \over 2}\right)^{r-1} 
\EE_{S_0}  \left| \EE_{S\setminus S_0} e^{\lambda f'} \right|. \endaligned \tag2.3.9$$
Combining (2.3.8) and (2.3.9), we conclude that 
$${\left| \EE_S\thinspace e^{\lambda  f''} - \EE_S\thinspace e^{\lambda f'} \right| \over  \left| \EE_S\thinspace e^{\lambda f'} \right|} \ \leq \ 
{2 |\lambda| \over \cos^{r-1} \left(\theta/2\right)}. $$
Recalling the bound for $|\lambda|$,  formula (2.3.1) for $\theta$ and using Part (1) of Lemma 2.2, we obtain 
$$\aligned {\left| \EE_S\thinspace e^{\lambda  f''} - \EE_S\thinspace e^{\lambda f'} \right| \over  \left| \EE_S\thinspace e^{\lambda f'} \right|} 
\ \leq \ &{2 \over  3c \sqrt{r-1}} \left(\cos {\pi \over 4 \sqrt{r-1}}\right)^{-(r-1)} \\ \leq \ &{2 \over  3c \sqrt{r-1}} {1 \over \cos (\pi/4)} 
= {2\sqrt{2} \over 3c \sqrt{(r-1)}}. \endaligned$$
It follows now that the angle, call it $\alpha$, between $\EE_S\thinspace e^{\lambda f'}$ and $\EE_S\thinspace e^{\lambda f''}$ is acute and that 
$$\sin \alpha \ \leq \  {2\sqrt{2} \over 3c \sqrt{(r-1)}}. $$
From (2.3.1) and Part (2) of Lemma 2.2, we have 
$$\sin {\theta \over c} \ = \ \sin {\pi \over 2 c \sqrt{r-1}} \ \geq \ {1 \over c \sqrt{r-1}} \sin {\pi \over 2} ={1 \over c \sqrt{r-1}} \ \geq \ \sin \alpha.$$
Hence the angle between $\EE_S\thinspace e^{\lambda f'}$ and $\EE_S\thinspace e^{\lambda f''}$ indeed does not exceed $\theta/c$ and $\EE_S\thinspace e^{\lambda f}$ rotates by not more than an angle of $\theta$ when the value of one of the coordinates $x_j$ with $j \notin S$ changes, while the others remain the same. This completes the induction step in proving Statement 2.3.2.

It remains to prove (1.1.2) assuming (1.1.1). Clearly, we have 
$$|f(x)| \ \leq \ m L \quad \text{for all} \quad x \in X$$
and the upper bound in (1.1.2) follows. To prove the lower bound, iterating (2.3.3), we get 
$$\left| \EE\thinspace e^{\lambda f} \right| \ \geq \ \left( \cos {\theta \over 2}\right)^n  \EE\thinspace \left| e^{\lambda f} \right| \ \geq \ 
e^{-|\lambda| m L } \left(\cos  {\pi \over 4 \sqrt{r-1}}\right)^n,$$
as required.
{\hfill \hfill \hfill} \qed

\head 3. Proof of Theorem 1.2 \endhead

We obtain Theorem 1.2 as a scaling limit of Theorem 1.1.

\subhead (3.1) Proof of Theorem 1.2 \endsubhead 
First, we consider the case of bounded functions $\phi_i$, so that the condition (1.2.1) is satisfied.

We are going to use Theorem 1.1. 
For an integer $N \geq 1$, we consider $nN$ probability spaces $X_1= \ldots = X_{nN}=\{-1, 1\}$ and their direct product 
$X=\{-1, 1\}^{nN}$, all endowed with the uniform probability measure. We write a point $x \in X$ as $x=\left(x_{jk}\right)$, where $j=1, \ldots, n$ and 
$k=1, \ldots, N$. For $i=1, \ldots, m$, we define functions $\psi_i: \{-1, 1\}^{nN} \longrightarrow {\Bbb C}$ by 
$$\psi_i(x)={\sqrt{N} \over 2} \phi_i\left( {x_{11} + \ldots + x_{1N} \over \sqrt{N}},\ldots, {x_{j1}+ \ldots + x_{jN} \over \sqrt{N}}, \ldots, {x_{n1} + \ldots + x_{nN} \over \sqrt{N}}\right)$$
and
$$g_N=\sum_{i=1}^m \psi_i.$$
Since the functions $\phi_i$ are 1-Lipschitz in the $\ell^1$ metric of ${\Bbb R}^n$, the functions $\psi_i$ are 1-Lipschitz in the Hamming metric of 
$X$. Moreover, each function $\psi_i$ depends on not more than $r_N=rN$ coordinates and at most $c$ functions $\psi_i$ depend on any particular coordinate $x_{jk}$. Finally, from (1.2.1), we conclude that 
$$\left| \psi_i(x) \right| \ \leq \ {\sqrt N \over 2} L \quad \text{for} \quad i=1, \ldots, m.$$
Therefore, from formula (1.1.2) of Theorem 1.1, we conclude that 
$$\left| \EE\thinspace e^{\lambda_N g_N} \right| \ \geq \ \exp\left\{-{1 \over 2} |\lambda_N| m \sqrt{N} L \right\} \left( \cos {\pi \over 4 \sqrt{rN-1}}\right)^{nN}.$$
provided 
$$|\lambda_N| \ \leq \ {1 \over 3 c \sqrt{r N}}.$$
Therefore,
$$\left| \EE\thinspace \exp\left\{ \lambda {2 \over \sqrt{N}} g_N \right\} \right| \ \geq \ \exp\left\{-|\lambda| m  L \right\} \left( \cos {\pi \over 4 \sqrt{rN-1}}\right)^{nN} \tag3.1.1$$
provided 
$$|\lambda| \ \leq \ {1 \over 6 c \sqrt{r}}.$$
By the Central Limit Theorem, as $N \longrightarrow \infty$, the random vector
$$\left({x_{11} + \ldots + x_{1N} \over \sqrt{N}},\ldots, {x_{j1}+ \ldots + x_{jN} \over \sqrt{N}}, \ldots, {x_{n1} + \ldots + x_{nN} \over \sqrt{N}}\right)$$
converges in distribution to the standard Gaussian measure in ${\Bbb R}^n$. Since the functions $\phi_i$ are continuous and bounded, 
we have 
$$\lim_{N \longrightarrow \infty}  \EE\thinspace \exp\left\{ \lambda {2 \over \sqrt{N}} g_N \right\} = \EE\thinspace e^{\lambda f},$$
where the expectation in the right hand side is taken with respect to the standard Gaussian measure in ${\Bbb R}^n$, see, for example, Section 7.2 of 
\cite{GS20}.  Since 
$$\lim_{N \longrightarrow \infty} \left( \cos {\pi \over 4 \sqrt{rN-1}}\right)^{nN}= \lim_{N \longrightarrow \infty} \left(1-{\pi^2 \over 32 (rN-1)}\right)^{nN}=
\exp\left\{ - {\pi^2 n\over 32 r} \right\},$$
from (3.1.1) we obtain the lower bound in the inequality (1.2.2). The upper bound in (1.2.2) is trivial.

It remains to consider the general case of not necessarily bounded functions $\phi_i$. Shifting, if necessary, 
$$\phi_i:= \phi_i - \phi_i(0) \quad \text{for} \quad i=1, \ldots, m,$$
without loss of generality we assume that 
$$\phi_i(0)=0 \quad \text{for} \quad i=1, \ldots, m.$$
Then 
$$\left|\phi_i\left(x_1, \ldots, x_n\right)\right| \ \leq \ \sum_{j=1}^n |x_j| \quad \text{for} \quad i=1, \ldots, m. \tag3.1.2$$
For $L>0$, we define the truncation of $\phi_i$ by
$$\phi_{i, L}(x)=\cases \phi_i(x) &\text{if}\quad  |\phi_i(x)| \leq L \\ L \phi_i(x)/  |\phi_i(x)|  &\text{if} \quad |\phi(x)| > L. \endcases$$
Then
$$\left| \phi_{i, L}(x) - \phi_{i, L} (y)\right| \ \leq \ \left| \phi_i(x) - \phi_i(y) \right|,$$
so $\phi_{i, L}$ satisfy the conditions of Theorem 1.2, and are bounded. Hence for 
$$f_L=\sum_{i=1}^m \phi_L$$
we have 
$$\EE\thinspace e^{\lambda f_L} \ne 0 \quad \text{provided} \quad |\lambda| \ \leq \ {1 \over 6 c \sqrt{r-1}}.$$
From (3.1.2) it follows that 
$$\lim_{L \longrightarrow +\infty} \EE\thinspace e^{\lambda f_L} = \EE\thinspace e^{\lambda f}$$ and that the convergence is uniform on 
any compact set of $\lambda$ in ${\Bbb C}$. Then by the Hurwitz Theorem, see, for example, Section 1 of \cite{Kr01}, we have two options: either $\EE\thinspace e^{\lambda f}=0$ 
for all $\lambda \in {\Bbb C}$ or $\EE\thinspace e^{\lambda f} \ne 0$ for all $\lambda$ in the open disc $|\lambda| < 1/6c \sqrt{r-1}$. Since for $\lambda=0$ we clearly 
have $\EE\thinspace e^{\lambda f} =1$, the first option is not realized.
{\hfill \hfill \hfill} \qed

\head 4. Optimality \endhead

Our goal is to show that the bound 
$$|\lambda| \ \leq \ {1 \over 3 c \sqrt{r-1}}$$
in Theorem 1.1 is optimal, up to a constant factor, replacing $1/3$. We do it in two steps. First, we show that for $c=1$ the dependence on $r$ is optimal, up to a constant factor, and then that the dependence on both $c$ and $r$ is optimal.

\subhead (4.1) Dependence on $r$ \endsubhead Let $X_1= \ldots = X_n =\{-1, 1\}$ and let $X=\{-1, 1\}^n$, all endowed with the uniform probability measure.
We intend to show that for all sufficiently large $n$ there is a function $\phi: X \longrightarrow {\Bbb C}$ which is 1-Lipschitz in the Hamming metric of $X$ and such that $\EE\thinspace e^{\lambda \phi} =0$ for some $\lambda$ satisfying 
$$|\lambda| \ \leq \ {\gamma \over \sqrt{n}},$$
where $\gamma >0$ is an absolute constant. 

We reverse engineer an example from our proof of Theorem 1.2 in Section 3.

The functions 
$$z \longmapsto \int_0^{+\infty} e^{z x} e^{-x^2/2} \ dx  \quad \text{and} \quad z \longmapsto  \int_{-\infty}^0 e^{zx} e^{-x^2/2} \ dx,$$
where we integrate over a real variable $x$,
are non-constant entire functions of $z \in {\Bbb C}$ and hence the range of each, by the Picard Theorem, is the whole plane ${\Bbb C}$, except perhaps one point. Therefore, there are two points $u, v \in {\Bbb C}$ such that 
$$\int_{-\infty}^0  e^{u x} e^{-x^2/2} \ dx + \int_0^{+\infty} e^{v x} e^{-x^2/2} \ dx=0.$$
Let us define $\psi: {\Bbb R} \longrightarrow {\Bbb C}$
$$\psi(x)=\cases v x &\text{if\ } x \geq 0 \\ u x &\text{if \ } x< 0. \endcases$$
Then 
$$\int_{-\infty}^{+\infty} e^{\psi(x)} e^{-x^2/2} \ dx =0. \tag4.1.1$$
It is easy to check that $\psi$ is Lipschitz, more precisely,
$$|\psi(x)-\psi(y)| \ \leq \ \tau |x-y| \quad \text{where} \quad \tau=\max\left\{ |u|, |v| \right\}.$$
For an $L > 0$, we consider the truncation 
$$\psi_L(x) =\cases \psi(x) & \text{if} \quad |\psi(x)| \leq L \\ L \psi(x)/|\psi(x)| &\text{if} \quad |\psi(x)| > L. \endcases$$
Then 
$$\left| \psi_L(x) - \psi_L(y)\right| \ \leq \ | \psi(x) - \psi(y)| \ \leq \ \tau |x-y|, \tag4.1.2$$
so $\psi_L$ is also Lipschitz and, in addition, bounded. 

Let us fix some $\rho > 1$.
We have 
$$\lim_{L \longrightarrow \infty} \int_{-\infty}^{+\infty} e^{z \psi_L(x)} e^{-x^2/2} \ dx = \int_{-\infty}^{+\infty} e^{z \psi(x)} e^{-x^2/2} \ dx \tag4.1.3$$
and from (4.1.2) it follows that the convergence is uniform in $z$ on all compact sets in ${\Bbb C}$. Since the right hand side of (4.1.3) for $z=0$ is equal to 1, while for $z=1$ is equal to 0 by (4.1.1), by the Hurwitz Theorem, see,
for example, Section 1 of \cite{Kr01}, we conclude that for all sufficiently large $L$ we must have 
$$\int_{-\infty}^{+\infty} e^{z \psi_L(x)} e^{-x^2/2} \ dx =0 \quad \text{for some} \quad z \quad \text{with} \quad |z| < \rho. \tag4.1.4$$

Let us pick some $L$ such that (4.1.4) holds. For an integer $n$, we consider $n$ probability spaces $X_1= \ldots = X_n= \{-1, 1\}^n$ and their product $X=\{-1, 1\}^n$, all endowed with the uniform probability measure. We define $\phi: X \longrightarrow {\Bbb C}$ by 
$$\phi(x) = {\sqrt{n} \over 2 \tau } \psi_L\left( {1 \over \sqrt{n}} \sum_{i=1}^n x_i \right) \quad \text{where} \quad x=\left(x_1, \ldots, x_n\right). \tag4.1.5$$
It follows from (4.1.2) that $\phi$ is 1-Lipschitz in the Hamming metric of $\{-1, 1\}^n$. Hence $\phi$ satisfies the conditions of Theorem 1.1 with 
$c=1$ and $r=n$.

Since by the CentraL Limit Theorem the normalized sum
$${1 \over \sqrt{n}}  \sum_{i=1}^n x_i $$
converges in distribution to the standard Gaussian measure as $n$ grows, and the function $\psi_L$ is continuous and bounded, we have 
$$\EE\thinspace \exp\left\{ z \psi_L \left({1 \over \sqrt{n}} \sum_{i=1}^n x_i \right)\right\} \longrightarrow {1 \over \sqrt{2 \pi}} \int_{-\infty}^{+\infty} e^{z \psi_L(x)} 
e^{-x^2/2} \ dx \tag4.1.6$$
and that the convergence in (4.1.6) is uniform in $z$ on all compact subsets of ${\Bbb C}$, see, for example, Section 7.2 of \cite{GS20}.
The right hand side of (4.1.6) is equal to 1 for $z=0$. Then by (4.1.4) and the Hurwitz Theorem, for all sufficiently large $n$ we must have 
$$\EE\thinspace \exp\left\{ z \psi_L \left({1 \over \sqrt{n}} \sum_{i=1}^n x_i \right)\right\} =0 \quad \text{for some} \quad z \quad \text{with} \quad |z|< \rho.$$
Therefore, from (4.1.5) we conclude that for all sufficiently large $n$, there is $\lambda$ satisfying 
$$|\lambda|  < {2 \tau \rho \over \sqrt{n}} \quad \text{and} \quad \EE\thinspace e^{\lambda \phi}=0,$$
which proves that the bound in Theorem 1.1 is indeed optimal in terms of $r$, up to an absolute constant factor. 

\subhead (4.2) Dependence on $r$ and $c$ \endsubhead
 To prove that the dependence on both $r$ and  $c$ is optimal, for an integer $k \geq 1$, we introduce
$$f_k(x)=\underbrace{\phi(x)+ \ldots + \phi(x)}_{\text{$k$ times}}, \tag4.2.1$$
where $\phi$ is defined by (4.1.5). While the parameter $r$ remains the same for all $f_k$ defined by (4.2.1), the parameter $c$ changes, $c=k$. Furthermore,
$$\EE\thinspace e^{\lambda f_k} = \EE\thinspace e^{(k\lambda) f},$$
and hence the zero-free region for $\lambda$ scales
$$|\lambda|=O\left({1 \over k \sqrt{r} }\right)=O\left({1 \over c \sqrt{r}}\right).$$
\subhead (4.3) Question: computational complexity \endsubhead It would be interesting to find out whether the dependence on $r$ is optimal from the computational complexity point of view, that is, whether for real-valued $f$ and $\lambda \gg 1/c\sqrt{r}$, the approximation of $\EE\thinspace e^{\lambda f}$
becomes computationally difficult, possibly conditioned on P $\ne$ NP or other commonly believed hypothesis. The argument with the ``cloning" of $f$ as in (4.2.1) shows that the dependence on $c$ is indeed optimal.

\head 5. Approximations \endhead 

\subhead (5.1) The goal \endsubhead 
Here we sketch how Theorems 1.1, respectively Theorem 1.2,  allow us to approximate $\EE\thinspace e^{\lambda f}$ provided 
$$|\lambda| \ \leq \ {1 -\delta \over 3 c \sqrt{r-1}}, \quad \text{respectively,} \quad |\lambda| \ \leq \ {1-\delta \over 6 c \sqrt{r-1}}\tag5.1.1$$
for some $0 < \delta < 1$, fixed in advance. The approach was used many times before, in particular in \cite{Ba17} in the context closest to ours.

The algorithm is based on the following result.
\proclaim{(5.2) Lemma} Let $p(z)$ be a univariate polynomial of degree $N$ in a complex variable $z$. Suppose that for some $\beta >1$, we have  $p(z)\ne 0$ for all $z$ satisfying 
$|z| < \beta$ and let us choose a branch of $g(z)=\ln p(z)$ in the disc $|z| < \beta$. For an integer $k \geq 1$, let $T_k(z)$ be the Taylor polynomial of $g(z)$ degree $k$ computed at $z=0$, that is, 
$$T_k(z)=g(0) + \sum_{s=1}^k {g^{(s)}(0) \over s!} z^s.$$
Then 
$$\left| g(z)-T_k(z)\right| \ \leq \ {N \over (k+1) (\beta-1) \beta^k} \quad \text{for all} \quad |z| \leq 1.$$
\endproclaim
\demo{Proof} This is Lemma 2.2.1 from \cite{Ba16}.
{\hfill \hfill \hfil} \qed
\enddemo

As follows from Lemma 5.2, to approximate $g(1)$ by $T_k(1)$ within an additive error $0 < \epsilon < 1$ and hence to approximate
$p(1)$ by $\exp\{T_k(1)\}$ within a relative error $0 < \epsilon < 1$, it suffices to choose 
$k=O_{\beta} \left(\ln N - \ln \epsilon\right)$. Moreover, the derivatives $g^{(s)}(0)$ can be computed from $p(0)$ and  the derivatives $p^{(s)}(0)$ for $s=1, \ldots, k$ in $O(k^2)$ time, see Section 2.2.2 of \cite{Ba16}.  

\subhead (5.3) Bounds \endsubhead
In the context of Theorem 1.1, let us pick an arbitrary $x_0 \in X$ and let $\gamma=f(x_0)$. In the context of Theorem 1.2, let $\gamma=f(0)$.
Since
$$\EE\thinspace e^{\lambda(f-\gamma)} =e^{-\lambda \gamma} \EE\thinspace e^{\lambda f},$$
without loss of generality, we assume that $f(x_0)=0$ (Theorem 1.1) or $f(0)=0$ (Theorem 1.2).

Then in the context of Theorem 1.1, respectively Theorem 1.2, we have 
$$\aligned &|f(x)| \ \leq \ nm \quad \text{for all} \quad x \in X , \quad \text{respectively,} \\ & |f(x_1, \ldots, x_n)| \ \leq \ m \sum_{j=1}^n |x_j| \quad \text{for all}
\quad (x_1, \ldots, x_n) \in {\Bbb R}^n. \endaligned \tag5.3.1$$

Given a complex $\lambda$ satisfying (5.1.1), we consider a function 
$$z \longmapsto \EE\thinspace e^{\lambda z f} = \EE\thinspace \exp\left\{ \lambda z \sum_{i=1}^m \phi_i \right\}.\tag5.3.2 $$ 
While (5.3.2) is not a polynomial, we can approximate it close enough by its Taylor polynomial $p_N(z)$ and then use Lemma 5.2 to approximate
$\ln p_N(z)$ by a polynomial of a low degree. To accomplish the first step, one can use, for example, the following standard estimate.

\proclaim{(5.4) Lemma} Let $\rho >0$ be a real number and let $N \geq 5\rho$ be an integer. Then
$$\left| e^z - \sum_{k=0}^N {z^k \over k!} \right| \ \leq \ e^{-2\rho} $$
for all $z \in {\Bbb C}$ such that $|z| \leq \rho$.
\endproclaim 
\demo{Proof} This is Lemma 3.2 in \cite{Ba17}. The proof is a standard estimate of the truncation error.
{\hfill \hfill \hfill} \qed
\enddemo

 Using (5.3.1), Lemma 5.4 and the bounds (1.1.2) and (1.2.2), for a given $0< \epsilon < 1$, one can compute 
$$N=\left((m+n) \ln {1 \over \epsilon} \right)^{O(1)},$$ such that the Taylor polynomial of (5.3.2),
$$p_N(z)=\sum_{s=0}^N {\lambda^s z^s \over s!} \EE f^s = \sum_{s=0}^N {\lambda^s z^s \over s!} \left( \sum_{i=1}^m \phi_i \right)^s$$  satisfies 
$$p_N(z) \ne 0 \quad \text{provided} \quad |z| \ \leq \ {1 \over1- \delta}$$ 
and $p_N(1)$ approximates $\EE\thinspace e^{\lambda f}$ within a relative error of $\epsilon/3$ (in the context of Theorem 1.2, we replace $\phi_i$ by their appropriate truncations). 

Using Lemma 5.2 with $\beta=(1-\delta)^{-1}$, 
we further approximate $p_N(1)$ within a relative error of $\epsilon/3$ using only $k=O_{\delta}\left(\ln (m+n)-\ln \epsilon)\right)$ first derivatives, 
$$p_N^{(s)}(0)=\lambda^s \EE \left(\sum_{i=1}^m \phi_i\right)^s \quad \text{for} \quad s=1, \ldots, k,$$
which in turn reduces to computing $m^{O_{\delta}\left( \ln (m+n)-\ln \epsilon\right)}$ different expectations
$$\EE \left(\phi_{i_1} \cdots \phi_{i_k}\right),$$
see, for example, \cite{Ba17} for a similar computation.

This results in a quasi-polynomial algorithm for approximating $\EE\thinspace e^{\lambda f}$. As we remarked, Patel and Regts \cite{PR17}, see also \cite{L+19}, developed methods for faster computation of the relevant derivatives of $g(z)=\ln p(z)$, which in some cases allows one to obtain a genuinely polynomial algorithm of 
$\bigl(m+n)/\epsilon\bigr)^{O(1)}$ complexity, provided the parameters $r$ and $c$ are fixed in advance.

\enddocument